\documentclass[letterpaper, 10 pt, conference]{ieeeconf}

\usepackage{graphicx} 
\usepackage{amsmath}
\usepackage{algorithm}
\usepackage{algpseudocode}
\usepackage{cite}

\title{\LARGE \bf
Tracking Point Vortices and Circulations via Advected Passive Particles: an Estimation Approach 
}

\author{Gil Marques$^{12}$, Marco Martins Afonso$^{1}$, Sílvio Gama$^{1}$}

\begin{document}

\maketitle
\thispagestyle{empty}
\pagestyle{empty}
\setcounter{footnote}{1}
\footnotetext{Gil Marques, Marco Martins Afonso and Sílvio Gama are with the Centro de Matemática da Universidade do Porto, Departamento de Matemática, Faculdade de Ciências, Universidade do Porto, Rua do Campo Alegre s/n, 4169-007, Porto, Portugal
        }
\setcounter{footnote}{2}
\footnotetext{{\tt\small gil.bj.marques@gmail.com}}

\begin{abstract}
We present a novel method for estimating the circulations and positions of point vortices using trajectory data of passive particles in the presence of Gaussian noise. The method comprises two algorithms: the first one calculates the vortex circulations, while the second one reconstructs the vortex trajectories. This reconstruction is done thanks to a hierarchy of optimization problems, involving the integration of systems of differential equations, over time sub-intervals all with the same amplitude defined by the autocorrelation function for the advected passive particles’ trajectories. Our findings indicate that accurately tracking the position of vortices and determining their circulations is achievable, even when passive particle trajectories are affected by noise.

\end{abstract}

\subsubsection*{Index Terms} -- Fluid flow systems, Optimization algorithms, Numerical algorithms

\section{Introduction}\label{sec:intro}

The objective of this work is to showcase a way to find the strength (or circulation) and trajectory of a system of point vortices by using the data of trajectories of passive particles. This is within the scope of the problems proposed by Protas \cite{Protas}, more specifically his Problem $4$, in which he considers the scenario of reconstructing some state of a vortex dynamical system by using incomplete or noisy data. 

The inverse problem of this, i.e. tracking passive particles advected by point vortices with known circulations is well documented in the literature. For example, the work \cite{boffetta1996trapping} uses Poincaré sections and canonical transformations to investigate the Hamiltonian dynamics of passive markers in the flow created by two point vortices. In \cite{huang2015detection}, the authors use data from coherent Lagrangian structures to determine when and where vortex formation occurs in flows around a 2D panel, and in \cite{qian2021tracking} the authors employ vortices to approximate the dynamics of an incompressible flow.
 
The study of point vortices and passive particles is interesting for understanding the global evolution of atmospheric and oceanic vortices~\cite{Meleshko1992,aref2001point, mokhov2020point}.

Section \ref{sec:ptvort} serves as an introduction to the formalism of point vortex dynamics. In Section \ref{sec:problem} we detail the problem we are tackling, presenting a method to solve it in Section \ref{sec:method}. We present the results obtained from testing the method on simulated corrupted data in Section \ref{sec:results}, and finally in Section \ref{sec:conc} we sum up some conclusions and ideas for future work.

\section{Point Vortices}\label{sec:ptvort}

In the vorticity formulation, the two-dimensional~(2D) incompressible Navier--Stokes equations can be written as
	
\begin{equation}
\frac{\partial\omega}{\partial t} + \left(\mathbf{u}\cdot\nabla\right)\omega = \nu\nabla^2\omega,	
\label{eq:vorticity2d}
\end{equation}
where $\omega$ is the vorticity, $\mathbf{u}$ is the velocity field and $\nu$ is the kinematic viscosity. In an inviscid flow ($\nu=0$), equation (\ref{eq:vorticity2d}), also known as Euler's equation, possesses the singular solutions 

\begin{equation}
\omega\left(x,y,t\right)=\sum_v\Gamma_v\,\delta\left(\left(x,y\right)-\left(x_v\left(t\right),y_v\left(t\right)\right)\right).
\end{equation}
The scalar $\Gamma_v$ denotes the circulation of the point vortex located in $\left(x_v\left(t\right),y_v\left(t\right)\right)$ and can be thought of as a measure of the rotation strength of the flow around the point vortex.


These singular solutions of the 2D incompressible Navier--Stokes equations \cite{Newton, Chorin, Batchelor, saffman1995vortex} were first studied by Helmholtz \cite{Helmholtz}, and further explored by Kelvin\cite{Kelvin1869} and Kirchhoff\cite{Kirchhoff1876}. They correspond to a scenario where the vorticity of a flow is concentrated in some well-defined points in space and enable a simpler description of the dynamics of such a system. Knowing the position and strength of each vortex in the system suffices to obtain the full velocity field at that instant, and thus knowing the dynamics of the vortices themselves in some time interval is enough to characterize the velocity field during that same time interval. 

While a single point vortex, if left alone, will stay in the position $\left(x_0,y_0\right)$ forever and its circulation $\Gamma$ will induce a time-independent velocity field in the plane, in a system of multiple point vortices the dynamics become more complex.

Consider a 2D inviscid flow in the complex plane that is described by $N_v$ point vortices located in the positions $z_{v}=x_{v} + iy_{v}\,(v=1,\dots,N_v\,; i^2=-1)$. The velocity field generated by each of these vortices will affect the remaining $N_v-1$ point vortices, making them move in the 2D plane, which causes the velocity field to change in time. Assuming there is no other object or force interfering in the flow, the motion of point vortices follows the differential equations

\begin{equation}
	\dot{z}_{v}^*=\frac{1}{2\pi i}\sum_{\substack{s=1\\ s\neq v}}^{N_v}\frac{\Gamma_{s}}{z_{v}-z_{s}}, \quad v=1,\dots,N_v\,,
\label{eq:mov_vortex}
\end{equation}
where $\Gamma_v$ is the circulation of vortex $v$ and $^*$ denotes the complex conjugate. 

However, when looking at real phenomena, it can be difficult to state precisely where the vortices are located and what is their circulation. A test particle, or a passive particle, is a particle that does not affect the flow and thus does not act on the vortices. Thus, they can be considered as point vortices with zero circulation. These passive particles, whose number we assume to be $N_p\,,$ can be easier to trace and the trajectory $z_p\left(\cdot\right)$ of any of them follows the equations of motion 
\begin{equation}
	\dot{z}_p^*=\frac{1}{2\pi i}\sum_{v=1}^{N_v}\frac{\Gamma_{v}}{z_{p}-z_{v}}, \quad p=1,\dots,N_p\,.
\label{eq:mov_part}
\end{equation}

\section{Problem}\label{sec:problem}
Assuming we have a way to track individual passive particles in a flow (for example, using drones, meteorological sensors, particle image/tracking velocimetry techniques), consider the problem of finding the circulations and positions of a system of $N_v$ point vortices, while only knowing the number of vortices in the system and the trajectory of $N_p$ passive particles advected by that system of vortices. 

In fact, since the trajectories of passive particles are sampled, what we have access to are discrete measurements of the particle trajectories. For simplicity, let us assume that these measurements are taken on $N_t$ equally spaced time instants $t_k = t_0 + hk$ in the time interval $\left[ t_0, t_f\right]$, where $h$ is the time increment.  

Thus, one way to approach the problem of reconstituting the vortices trajectories is to use the information on the particle trajectories to trace back the vortices that originated them, by trying to reconstitute the measured particle trajectories using the knowledge of the equations of motion for both the particles and the vortices. This can be summed up in the following problem, where $\widetilde{z}_{p}\left(t\right)$ is the position of particle $p$ at time $t$ as obtained from solving (\ref{eq:mov_part}): 

\begin{equation}
\left\{
\begin{aligned}
    \text{Minimize} &\quad \sum_{k=0}^{N_t-1}\sum_{p=1}^{N_p} || \widetilde{z}_p\left(t_k\right) - z_p\left(t_k\right) ||^2 \\
    \text{subject to}& \\
        &\dot{\widetilde{z}}_{p}^*=\frac{1}{2\pi i}\sum_{v=1}^{N_v}\frac{\Gamma_{v}}{\widetilde{z}_{p}-z_{v}}, \quad p=1,\dots,N_p \\
        &\dot{z}_{v}^*=\frac{1}{2\pi i}\sum_{\substack{s=1\\ s\neq v}}^{N_v}\frac{\Gamma_{s}}{z_{v}-z_{s}}, \quad v=1,\dots,N_v \\
        &\widetilde{z}_{p}\left(t_0\right) = z_{p}\left(t_0\right) \\
        &z_v\left(t_0\right) \quad \text{free} \\
        &\Gamma_v \qquad\!\quad \text{free}
\end{aligned}
\right.,
\label{eq:problem}
\end{equation}
where $||\cdot||$ is the usual Euclidean norm.

As a consequence to solving this problem, we will obtain the information about the vortices circulations and trajectories.

\section{Method}\label{sec:method}

Notice that first and foremost, since the particle trajectory data is prone to unknown measurement errors, we actually do not measure the real $z_p\left(t_k\right)$, but a quantity $z_p\left(t_k\right) + \Delta z_p\left(t_k\right)$, where $\Delta z_p\left(t\right)$ is the error associated with that measurement. By using a smoothing algorithm to diminish the impact of measurement errors in the data, we can then obtain estimations for the particle trajectories. From here onward, $z_p$ will be used to denote the filtered data in order to not overload the notation. 

From the filtered trajectories one can use a finite-differences method to estimate the velocity of the passive particles $\dot{z}_p\left(t\right)$. However, if velocity data are available, one can smooth such data directly instead of estimating them from the filtered trajectory data, which is something that would introduce further errors in the velocity estimates.

We can take advantage of the fact that $\Gamma_v$ is constant in time to divide the problem into two sub-problems. First, we will find the circulations of the vortices, and afterwards use them to find better estimates for the vortices trajectories.

\subsection{Estimating the circulations}

For any given time $t$ for which we have data for, the unknowns in (\ref{eq:mov_part}) are $\Gamma_v$ and $z_v\left(t\right)$. These amount, in reality, to $3N_v$ unknowns and thus we need to know the trajectories of at least
$3N_v/2$
particles (each brings two pieces of information) in order to build a system of equations that is not underdetermined. After building such a system, one can use a nonlinear solver method to find a solution of the system. This will require an initial guess for both the circulations and positions of the vortices that should be provided based on analysis of the data: either the movement of the particles in the flow or actual visual imagery of the environment should hint towards what range of values one can expect for these quantities. 

Doing so for time $t_k$ will return estimates for the circulations $\Gamma_v$ and for $z_v\left(t_k\right)$ --- the positions of the vortices at time $t_k$. We will identify these estimates as $\widetilde{\Gamma}_{v_k}$ and $\widetilde{z}_{v_k}$. The index $k$ on the circulations is used just to identify that such a value was estimated using the data at time $t_k\,,$ as the circulation is in fact independent of time.

One can then use $\widetilde{z}_{v_k}$ as an initial condition for (\ref{eq:mov_vortex}) in order to find an initial guess for $z_v\left(t_{k+1}\right)$ and build the same nonlinear system of equations as before for time $t_{k+1}$, solve it and obtain estimates for $z_v\left(t_{k+1}\right)$ and a new estimate for $\Gamma_v$: $\widetilde{\Gamma}_{v_{k+1}}$. In order to explore the parameter space further, one can update the initial guess for the circulations as $\widetilde{\Gamma}_{v_k}\left( 1 + \delta_v\right)$, where $\delta_v$ is a uniform random variable that ranges from $-\epsilon$ to $\epsilon$ ($\epsilon$ is some small positive number).

Repeating this procedure from $t_0$ to $t_f$ will result in a set of values $\left\{\widetilde{\Gamma}_{v_k}:k=1,\dots,N_t\right\}$, all of which are estimates for the vortex circulations. By doing an outlier analysis one can remove any abnormal value and find estimates for the vortex circulations using the average or median of the leftover values in order to obtain $\widetilde{\Gamma}_v$, a more robust estimate for these quantities. One can even repeat this whole procedure using these more robust estimates as the initial guess in order to refine the results. Algorithm \ref{alg:circ} sums up this procedure.

\begin{algorithm}[b]
\caption{Finding the circulations}\label{alg:circ}
\begin{algorithmic}[1]
\State Set a value for $\epsilon$
\State Estimate the velocities of passive particles
\State Provide an initial guess for the circulations and positions of the vortices at $t=t_0$
\State $k \gets 1$
\While{ $k \leq N_t$ }
\State Use a nonlinear solver to find $\widetilde{\Gamma}_{v_k}$ and $\widetilde{z}_{v_k}$, estimates for $\Gamma_v$ and $z_v\left(t_k\right)$ with the appropriate initial conditions
\State Use $\widetilde{z}_{v_k}$ and (\ref{eq:mov_vortex}) to find an estimate for $z_v\left(t_{k+1}\right)$
\State Generate random numbers $\delta_v$ uniformly in the interval $\left[-\epsilon,\epsilon\right]$
\State Use $\widetilde{\Gamma}_{v_k}\left( 1 + \delta_v\right)$ as an updated initial guess for the circulations
\State $k \gets k+1$
\EndWhile
\State Run an outlier analysis over the values $\widetilde{\Gamma}_{v_k}$
\State Obtain a more robust estimate of the circulations using the average or median of the leftover $\widetilde{\Gamma}_{v_k}$ values
\end{algorithmic}
\end{algorithm}

\subsection{Estimating the vortex trajectories}

Having obtained the circulations, we have now simplified the problem (\ref{eq:problem}) as $\Gamma_v$ is now fixed at $\widetilde{\Gamma}_v$. 

It is known that, in general, vortex systems exhibit chaotic dynamics for $N_v\geq 4$ \cite{Babiano,Newton}. A practical consequence of this fact is that, when solving equations such as (\ref{eq:problem}) or (\ref{eq:mov_part}), even a small error in the initial condition can be highly amplified given enough time and the computed trajectories can be completely different from the real ones. Our approach to deal with this problem is to partition the trajectory and work on the reconstitution of each part separately and sequentially. However, one needs to first find a criteria for building such a partition which should take into account what is known about the chaoticity of the system.

 One way to measure the degree of chaos in such a system are Lyapunov exponents\cite{ruelle1989chaotic,Chaos}. In particular, the largest Lyapunov exponent measures the (exponential) rate at which two nearby trajectories separate from each other, and thus its inverse gives us the time needed for two nearby trajectories to separate for a factor of $e$. This also tells us that the system has some degree of memory but, after a certain time, the evolution is mostly uncorrelated with the states past a certain point in time. Lyapunov exponents for the system of vortices and for passive particles with a chaotic movement are different. Nevertheless, in a chaotic system, they seem to be within the same order of magnitude \cite{controlo2020}. Thus, since we have no information about the vortices themselves, we need to use the passive particles to get an estimate of such quantities. However, Lyapunov exponents computation is slow and heavy and, in a scenario where one wants to follow vortices in real time, it is not a quantity suitable to use. Autocorrelation functions, nonetheless, are fast to compute and can be used to find the time after which the state of the system has little to no memory of its initial state. 

 By computing autocorrelation functions for the passive particles' trajectories, we can find the time $\tau$ it takes for the system to lose most of its memory from the initial state. We do this by setting an arbitrary cutoff $0<\alpha<1$ and computing $\tau : \forall p \in \left\{1,\dots,N_p\right\} ,\,|\,\rho_p\left(t>\tau\right)|\leq\alpha\,,$ where $\rho_p$ denotes the autocorrelation function for the trajectory of the particle~$p$. 

 Thus, in order to reconstitute the full vortex trajectories we must reconstitute $n=\Bigr\lceil\frac{t_f-t_0}{\tau}\Bigr\rceil$ partitions of the form $P_j=\left[t_0+j\tau, t_0+\left(j+1\right)\tau\right],\, j=0,1,\dots, n-2\,,$ and $P_{n-1} = \left[t_0+\left(n-1\right)\tau,t_f\right]$. Each of these constitutes an optimization problem akin to (\ref{eq:problem}) that can be written as

 \begin{equation}
\left\{
\begin{aligned}
    &\text{Minimize} \quad \sum_{k=j\tau/h}^{\left(j+1\right)\tau/h}\sum_{p=1}^{N_p} || \widetilde{z}_p\left(t_k\right) - z_p\left(t_k\right) ||^2 \\
    &\text{subject to} \\
        &\dot{\widetilde{z}}_{p}^*=\frac{1}{2\pi i}\sum_{v=1}^{N_v}\frac{\widetilde{\Gamma}_{v}}{\widetilde{z}_{p}-z_{v}}, \quad p=1,\dots,N_p \\
        &\dot{z}_{v}^*=\frac{1}{2\pi i}\sum_{\substack{s=1\\ s\neq v}}^{N_v}\frac{\widetilde{\Gamma}_{s}}{z_{v}-z_{s}}, \quad v=1,\dots,N_v \\
        &\widetilde{z}_{p}\left(t_0+ j\tau\right) = z_{p}\left(t_0 + j\tau\right) \\
        &z_v\left(t_0 + j\tau\right) \quad \text{free} \\
        & t\in P_j
\end{aligned}
\right.
\label{eq:subproblem}
\end{equation}
for $j=0,1,\dots,n-2\,,$ and

 \begin{equation}
\left\{
\begin{aligned}
    &\text{Minimize} \, \sum_{k=\left(n-1\right)\tau/h}^{N_t-1}\sum_{p=1}^{N_p} || \widetilde{z}_p\left(t_k\right) - z_p\left(t_k\right) ||^2 \\
    &\text{subject to} \\
        &\dot{\widetilde{z}}_{p}^*=\frac{1}{2\pi i}\sum_{v=1}^{N_v}\frac{\widetilde{\Gamma}_{v}}{\widetilde{z}_{p}-z_{v}}, \quad p=1,\dots,N_p \\
        &\dot{z}_{v}^*=\frac{1}{2\pi i}\sum_{\substack{s=1\\ s\neq v}}^{N_v}\frac{\widetilde{\Gamma}_{s}}{z_{v}-z_{s}}, \quad v=1,\dots,N_v \\
        &\widetilde{z}_{p}\left(t_0+ \left(n-1\right)\tau\right) = z_{p}\left(t_0 + \left(n-1\right)\tau\right) \\
        &z_v\left(t_0 + \left(n-1\right)\tau\right) \quad \text{free} \\
        &t\in P_{n-1}
\end{aligned}
\right.
\label{eq:subproblem2}
\end{equation}

These $n$ problems can thus be solved sequentially by a nonlinear solver, using $z_v\left(t_0 + \left(j+1\right)\tau\right)$ obtained in the problem associated to the interval $P_j$ as an initial guess for the problem associated with the interval $P_{j+1}$. Algorithm \ref{alg:traj} summarizes this full process.

\begin{algorithm}[b]
\caption{Finding the vortex trajectories}\label{alg:traj}
\begin{algorithmic}[1]
\State Set $\alpha$
\State $p \gets 1$
\While{ $p \leq N_p$ }
\State Compute $\rho_p\left(t\right)$, the autocorrelation function of $z_p\left(t\right)$
\State Find $\tau_p : \,|\,\rho_p\left(t>\tau_p\right)|\leq\alpha$
\State $p \gets p + 1$
\EndWhile
\State $\tau \gets \min_{1\leq p\leq N_p}\left\{\tau_p\right\}$
\State $n \gets \Bigr\lceil\frac{t_f-t_0}{\tau}\Bigr\rceil$
\State $j \gets 1$
\While{ $j \leq n$ }
\State Use a nonlinear solver to solve problem (\ref{eq:subproblem}) on the time interval $P_j$
\State $j \gets j+1$
\EndWhile
\State Use a nonlinear solver to solve problem (\ref{eq:subproblem2})
\end{algorithmic}
\end{algorithm}

\section{Results}\label{sec:results}

We simulated a system comprised of $N_v=4$ vortices with circulations $\Gamma_v=v,\,v=1,2,3,4$ and $N_p=20$ passive particles. The vortices were initially placed in the position $z_1\left(t=0\right) = 2$, $z_2\left(t=0\right) = -1-i$, $z_3\left(t=0\right) = \frac{1}{2}+\frac{i}{2}$, $z_4\left(t=0\right) = -2+3i$ and the passive particles were placed in random positions in the background between the vortices, and we let the system evolve for a time period of $10^5$ units of time using a time step of $10^{-2}$. The vortex circulations were all chosen to have the same sign in order to confine the movements of the particles and vortices to a certain area, which facilitates the simulation and analysis. If there are vortices with opposite signs in the system, pairs of vortices can be formed and travel long distances away from the remaining vortices. The simulation ran for a long time in order to be possible for us to compute a good estimate of the maximum Lyapunov exponent of the vortex system, for which we found $\sim 0.028$. However, since the problem of recovering the vortex trajectories is computationally heavy, we restrained to using a time interval of $10^3$ units of time for this problem.

We then introduced Gaussian white noise $W\sim \mathcal{N}\left(\mu=0,\sigma^2=0.01^2\right)$ in the particles' trajectory data, in order to simulate what actual real data could look like. After doing so, the data were too corrupted to be directly used, so we first smoothed them using a Gaussian filter. We then used the smoothed trajectories to compute particle velocities using a fourth-order finite-differences scheme and ran Algorithm \ref{alg:circ} with $\epsilon=0.1$. Integration was made using the Runge-–Kutta-–Fehlberg method with a step size of $10^{-2}$, and the nonlinear equations were solved by 
MATLAB$^{TM}$ routine \textit{lsqcurvefit}.

We were able to recover the vortices' circulations with a relative error of $\sim 0.1\%\,$ (if we use the uncorrupted data, the relative error decreases to $\sim 10^{-7}$). The obtained values are summarized in Table \ref{tb:circ}, as well as the relative errors for each individual circulation. 

\begin{figure*}[t]
	\centering
    \vspace{5mm}
	\includegraphics[width=0.45\linewidth]{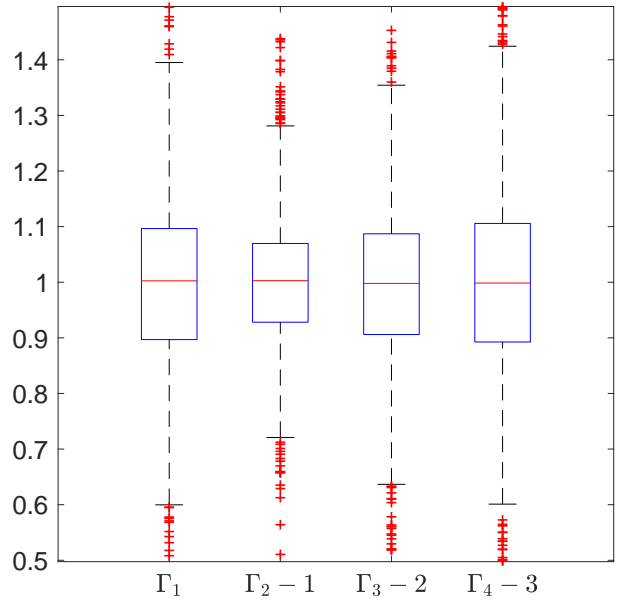}
	\caption[Boxplot]{Boxplot for $\widetilde{\Gamma}_v - (v-1)$. The values are shifted down by $v-1$ for easier visualization purposes. We can see that the distribution of obtained values is centered around the true value of each circulation.}
	\label{fig:boxplot}
\end{figure*}

Fig. \ref{fig:boxplot} shows the boxplots for the values of $\widetilde{\Gamma}_{v_k}$ obtained by the algorithm. The values in the~$y$ axis were shifted down by $v-1$ for easier visualization purposes.

\begin{table*}[t]
    \vspace{5mm}
\caption{True values of the circulations and the values recovered by Algorithm \ref{alg:circ}. The circulations were recovered with a relative error of $\sim0.1\%$.}
\label{tb:circ}
\begin{center}
\begin{tabular}{c|cccc}
Vortex $v$ & $1$ & $2$ & $3$ & $4$\\
\hline
$z_v\left(0\right)$ & $2$ & $-1-i$ & $\frac{1}{2}+\frac{i}{2}$ & $-2+3i$ \\
$\Gamma_v$ & $1$ & $2$ & $3$ & $4$\\
$\widetilde{\Gamma}_v$ & $0.99625$ & $2.00132$ & $2.99746$ & $4.00258$ \\
$\frac{|\widetilde{\Gamma}_v - \Gamma_v|}{|\Gamma_v|}$ & 
$3.7\times10^{-3}$ & $6.6\times10^{-4}$ & $8.5\times10^{-4}$ & $6.5\times10^{-4}$
\\
\end{tabular}
\end{center}
\end{table*}

Using the obtained values of $\widetilde{\Gamma}_v$, we then proceeded to recover the vortex trajectories using Algorithm \ref{alg:traj}. We set $\alpha=0.2$, obtaining $\tau=39.89\,.$ This is in accordance with the value of the Lyapunov exponent we obtained in the initial simulation ($\sim0.028$), since $\tau^{-1}\approx0.025\,,$ which is of the same order of magnitude as the Lyapunov exponent. Once again, integration was made using the Runge–-Kutta–-Fehlberg method with a step size of $10^{-2}$ and the nonlinear problems were solved by MATLAB$^{TM}$ using the routine \textit{fmincon} with the \textit{sqp} algorithm. We were able to recover the trajectories of the vortex system with a total relative error of  $\sim 0.44\%$. (if we use the uncorrupted data, the relative error decreases to $\sim 8\times 10^{-4}$ ). 

\begin{figure*}[t]
	\centering
     \vspace{5mm}

	\includegraphics[width=0.6\linewidth]{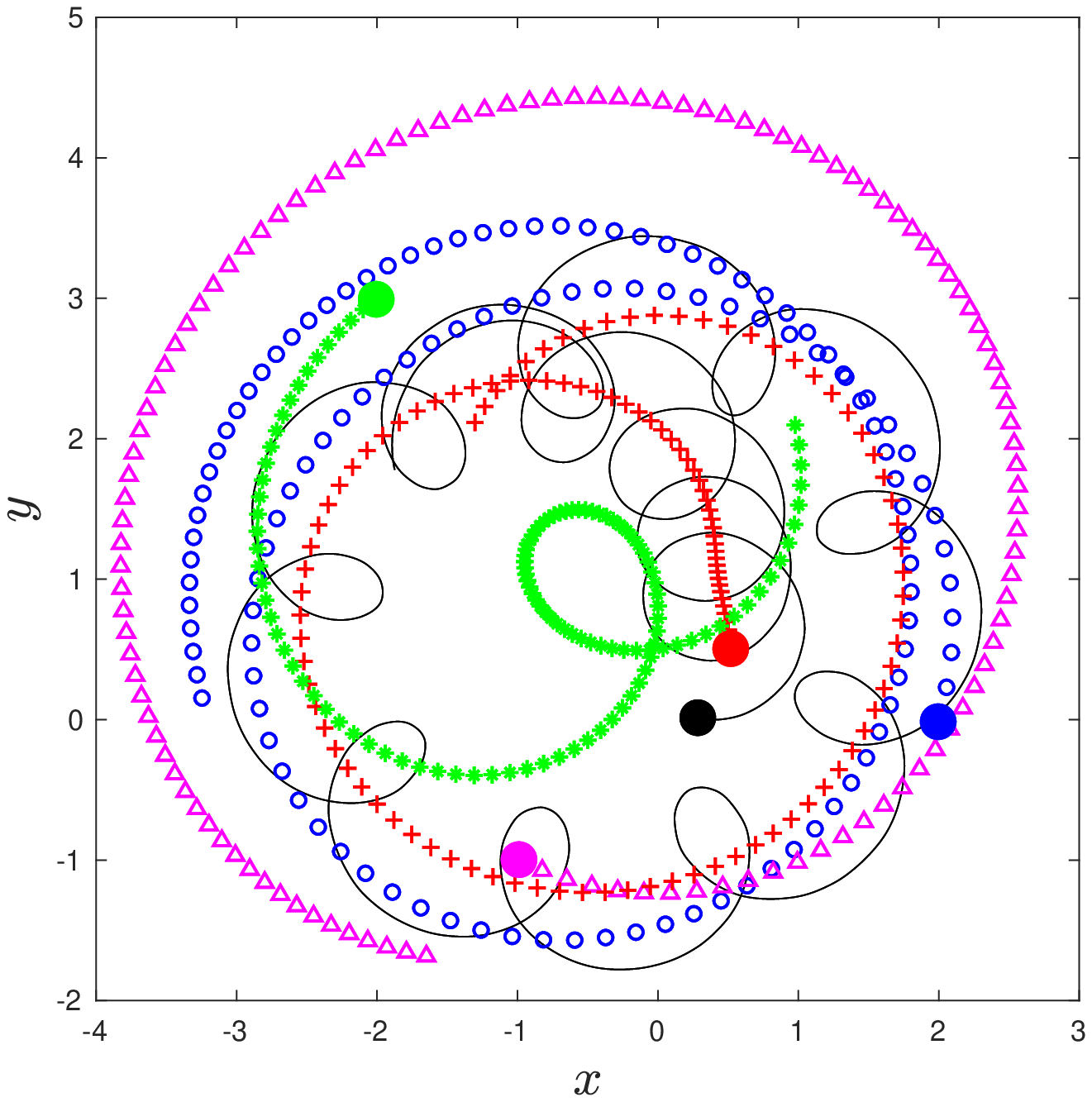}
	\caption[]{A portion of the recovered trajectories of the $4$ vortices (blue dots, pink triangles, red crosses and green asterisks identify vortices $1$, $2$, $3$ and $4$, respectively) and an example of a smoothed trajectory for a passive particle in a full black line. The color-filled dots mark the initial point of each trajectory.}
	\label{fig:Traj}
\end{figure*}

Fig. \ref{fig:Traj} shows a portion of the recovered trajectories of the $4$ vortices with different colored markers and the (smoothed) trajectory of a passive particle in a full black line. The color-filled dots mark the initial point of each trajectory. 

\begin{figure}[t]
	\centering    
 \vspace{5mm}

	\includegraphics[width=\linewidth]{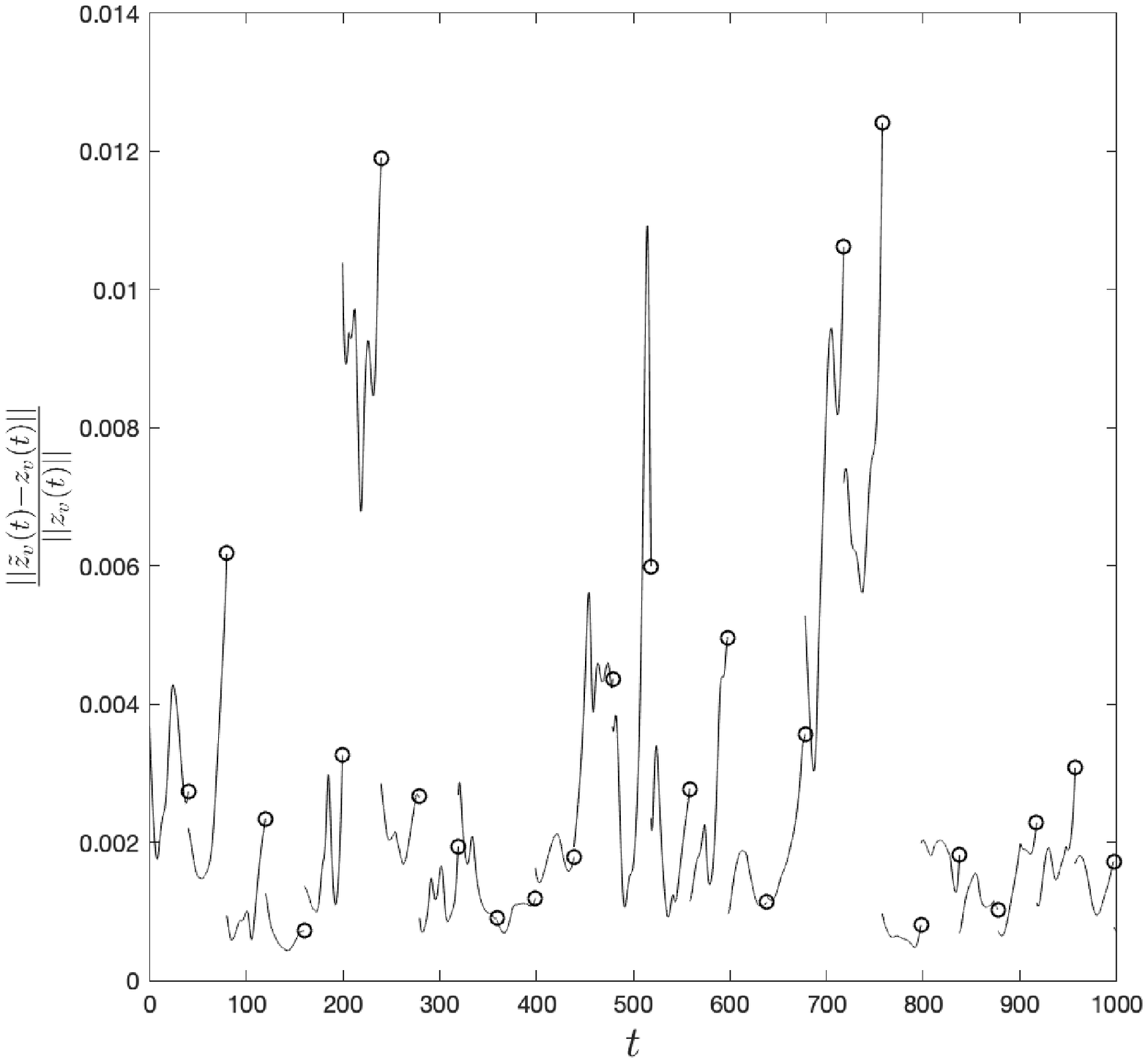}
	\caption[]{Relative error for the position of the vortex system at each time step. The circles identify the discontinuities that arise from the partition of the problem into various sub-problems. In general the relative error decreases every time we change from one partition to the following one.}
	\label{fig:TrajErr}
\end{figure}

Fig. \ref{fig:TrajErr} illustrates how the partitioning of the problem helps with controlling the error of the recovered trajectories. Due to the chaoticity of the system, errors can grow large and we can see that, in general, by integrating partitions of the trajectory with a time length similar to the inverse of the Lyapunov exponent of the vortex system will help to prevent exponential growth of the errors.

\section{Conclusion}\label{sec:conc}
We presented the basic ideas on how to try to recover vortex dynamics from data on the trajectories of advected passive particles. This is a problem that, as far as we are aware, has not been addressed before and can have relevant interest since it is easier to track the advected particles than tracking and measuring the strength of vortices themselves. 

The method presented is a first study and we believe it can be improved largely. For example, using optimization techniques and algorithms specifically tailored for the problem and its equations, instead of using generic MATLAB$^{TM}$ routines, will likely increase accuracy and speed up the computations. 

Nevertheless, from the results of this study we see that it is possible to effectively recover the full vortex dynamics with a good degree of accuracy, even in the presence of noisy data. This opens up the possibility of studying more complex problems, for instance, predicting positions of vortices and/or particles at future times, which enables the study of finding controls to move passive particles to a desired position in a given time.

\section*{Acknowledgements}
This work was supported by (i) CMUP, member of LASI, which is financed by national funds through FCT – Fundação para a Ciência e a Tecnologia, I.P., under the project with reference UIDB/00144/2020, and (ii) project SNAP NORTE-01-0145-FEDER-000085, co-financed by the European Regional Development Fund (ERDF) through the North Portugal Regional Operational Programme (NORTE2020) under Portugal 2020 Partnership Agreement. GM thanks grant ref. PD/BD/150537/2019 through FCT.



	\bibliographystyle{ieeetr}
	\bibliography{biblio}

\begin{thebibliography}{10}

\bibitem{Protas}
B.~Protas, ``Vortex dynamics models in flow control problems,'' {\em
  Nonlinearity}, vol.~21, no.~9, pp.~R203--R250, 2008.

\bibitem{boffetta1996trapping}
G.~Boffetta, A.~Celani, and P.~Franzese, ``Trapping of passive tracers in a
  point vortex system,'' {\em Journal of Physics A: Mathematical and General},
  vol.~29, no.~14, p.~3749, 1996.

\bibitem{huang2015detection}
Y.~Huang and M.~A. Green, ``Detection and tracking of vortex phenomena using
  {L}agrangian coherent structures,'' {\em Experiments in Fluids}, vol.~56,
  no.~7, p.~147, 2015.

\bibitem{qian2021tracking}
Z.~Qian, Y.~Qiu, and Y.~Zhang, ``Tracking the vortex motion by using {B}rownian
  fluid particles,'' {\em Physics of Fluids}, vol.~33, no.~10, p.~105113, 2021.

\bibitem{Meleshko1992}
V.~V. Meleshko, M.~Y. Konstantinov, A.~A. Gurzhi, and T.~P. Konovaljuk,
  ``Advection of a vortex pair atmosphere in a velocity field of point
  vortices,'' {\em Physics of Fluids A: Fluid Dynamics}, vol.~4, no.~12,
  pp.~2779--2797, 1992.

\bibitem{aref2001point}
H.~Aref and M.~A. Stremler, ``Point vortex models and the dynamics of strong
  vortices in the atmosphere and oceans,'' in {\em Fluid Mechanics and the
  Environment: Dynamical Approaches: A Collection of Research Papers Written in
  Commemoration of the 60th Birthday of Sidney Leibovich}, pp.~1--17, Springer,
  2001.

\bibitem{mokhov2020point}
I.~I. Mokhov, S.~G. Chefranov, and A.~G. Chefranov, ``Point vortices dynamics
  on a rotating sphere and modeling of global atmospheric vortices
  interaction,'' {\em Physics of Fluids}, vol.~32, no.~10, p.~106605, 2020.

\bibitem{Newton}
P.~Newton, {\em The N-Vortex Problem: Analytical Techniques}.
\newblock Appl. Math. Sci., Springer New York, 2001.

\bibitem{Chorin}
A.~Chorin, {\em Vorticity and Turbulence}.
\newblock Appl. Math. Sci., Springer, 1994.

\bibitem{Batchelor}
G.~K. Batchelor, {\em An Introduction to Fluid Dynamics}.
\newblock Cambridge Mathematical Library, Cambridge University Press, 2000.

\bibitem{saffman1995vortex}
P.~G. Saffman, {\em Vortex dynamics}.
\newblock Cambridge University Press, 1995.

\bibitem{Helmholtz}
H.~Helmholtz, ``Über {I}ntegrale der hydrodynamischen {G}leichungen, welche
  den {W}irbelbewegungen entsprechen.,'' vol.~1858, no.~55, pp.~25--55, 1858.

\bibitem{Kelvin1869}
W.~Thomson (Lord~Kelvin), ``On vortex motion.,'' {\em Trans. R. Soc. Edin},
  vol.~25, p.~217–260, 1869.

\bibitem{Kirchhoff1876}
G.~R. Kirchhoff, ``Vorlesungenb\"er mathematische {P}hysik.,'' {\em Mechanik},
  1876.

\bibitem{Babiano}
A.~Babiano, G.~Boffetta, A.~Provenzale, and A.~Vulpiani, ``Chaotic advection in
  point vortex models and two‐dimensional turbulence,'' {\em Phys. of
  Fluids}, vol.~6, no.~7, pp.~2465--2474, 1994.

\bibitem{ruelle1989chaotic}
D.~Ruelle and S.~Isola, {\em Chaotic evolution and strange attractors}, vol.~1.
\newblock Cambridge University Press, 1989.

\bibitem{Chaos}
K.~T. Alligood, T.~D. Sauer, and J.~A. Yorke, {\em Chaos An Introduction to
  Dynamical Systems}.
\newblock Springer-Verlag, New York, 1997.

\bibitem{controlo2020}
G.~Marques, M.~J. Rodrigues, and S.~Gama, ``Passive particle dynamics in
  viscous vortex flow,'' in {\em CONTROLO 2020} (J.~A. Gon{\c{c}}alves,
  M.~Braz-C{\'e}sar, and J.~P. Coelho, eds.), (Cham), pp.~352--362, Springer
  International Publishing, 2021.

\end{thebibliography}


\end{document}